\documentclass[letterpaper,10pt,nonacm]{article} 

\usepackage{fancyhdr}

\usepackage{cite}

\usepackage[sort, numbers]{natbib}

\usepackage{tikz}
\usetikzlibrary{shapes.geometric, arrows}

\usepackage{multirow}

\usepackage{subcaption}
\tikzstyle{startstop} = [rectangle, minimum width =3cm, minimum height = 1cm, text centered, draw=black, fill=white]
\tikzstyle{decision} = [diamond, minimum width = 3cm, minimum height = 1cm, text centered, text width = 2cm, draw=black, fill=white]
\tikzstyle{arrow} = [thick,->,>=stealth]

\usepackage{osameet3} 

\usepackage{algorithm}

\providecommand{\keywords}[1]
{
	\small	
	\textbf{\textit{Keywords---}} #1
}

\usepackage{amsmath,amssymb}
\usepackage[colorlinks=true,bookmarks=false,citecolor=blue,urlcolor=blue]{hyperref} 

\numberwithin{equation}{section}
\numberwithin{algorithm}{section}
\numberwithin{figure}{section}
\numberwithin{table}{section}

\usepackage{todonotes}
\usepackage{nicefrac}

\DeclareSymbolFont{epsilon}{OML}{ntxmi}{m}{it}
\DeclareMathSymbol{\epsilon}{\mathord}{epsilon}{"0F}

\begin{document}

\title{Efficient Geometry Representation Strategies for the Shape Optimization of Profile Extrusion Dies}

\author{Jana Sasse, Maximilian Esser, Markus Mügge, Stefan Turek}
\address{Institute for Applied Mathematics, LSIII, TU Dortmund University, D-44227 Dortmund Germany}
\email{jana.sasse@math.tu-dortmund.de}

\vspace{4pt}
\begin{minipage}{0.14\textwidth}
	\includegraphics[width=\textwidth]{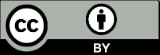}
\end{minipage}\hspace*{1em}
\begin{minipage}{0.84\textwidth}
	{\small This work is licensed under a Creative Commons Attribution 4.0 International\newline license (see \url{https://creativecommons.org/licenses/by/4.0/legalcode.de})}
\end{minipage}

\begin{abstract}
The design of profile extrusion dies remains a challenging task due to the complex rheological behavior of polymer melts and the geometric intricacies of flow channels.
Traditional manual optimization approaches, which rely heavily on human experience, are inefficient and often employ unvalidated heuristics.
To address these challenges, we present a deterministic and explainable framework for automatic die design based on adjoint-based shape optimization.
This approach enables the computation of sensitivities that directly indicate beneficial modifications to the flow channel geometry.
A major difficulty in such optimization processes lies in generating boundary-conforming meshes that evolve consistently with changing geometries.
To overcome this issue, we employ non-boundary conforming geometry representation methods that eliminate the need for an explicit surface representation along physical boundaries.
A dedicated reconstruction technique is developed to recover accurate sensitivity information at the virtual interface between fluid and solid regions.
The proposed algorithm is demonstrated on 3D geometries with varying complexity, including realistic extrusion die flow channels.
Several objective functionals relevant to industrial applications, such as flow balance at the outflow, are considered.
The results highlight significant improvements in performance metrics while maintaining numerical robustness.
This work showcases the potential of adjoint-based techniques for automated die design in a domain still largely governed by manual trial-and-error procedures, establishing a foundation for data-efficient, sustainable manufacturing workflows using computational rheology.
\end{abstract}
\hspace{10pt}

\vspace{16pt}


\section{Introduction}

The design of the flow channel in profile extrusion dies remains challenging. The optimal flow channel geometry creates a uniform outflow velocity profile. However, this is highly dependent on the operating point, i.e., the throughput, processing temperature, and thermo-rheological properties of the plastic melt. As a result, extrusion die design still relies heavily on manual design iterations.

The last years have seen an emergence of computational methods for extrusion die design.
However, these methods face various numerical challenges.
For complex 3D geometries, the modeling of the computational domain alone is a non-trivial task.	
In addition, the extruded plastic melt exhibits a complex rheological behavior, since it is shear-thinning, viscoelastic and temperature-dependent. 
Furthermore, the computational cost of established optimization methods~\cite{GasparCunha2022} is also significant due to the associated number of solver calls.
Adjoint-based sensitivities are a promising alternative at low computational cost~\cite{Othmer2008, Bletsos2023, Sasse2025}.

The goal of the research presented in this paper is to identify the most efficient method for adjoint-based optimization of complex profile extrusion dies.
This includes the most efficient geometry representation strategy, but also the best strategy for the sensitivity-based optimization algorithm.

\section{Numerical Modeling}

\subsection{Geometry Representation Strategies}\label{sec:Geometry}

The most common method employs a body-fitted mesh that describes the entire fluid domain.
While this method has the advantage that it operates on the "exact representation" of the flow channel geometry, it also comes with its drawbacks.
For complex flow channel geometries, the computational cost of the mesh generation becomes prohibitively expensive.
In addition, large deformations of the flow channel during the optimization process often require re-meshing during runtime.
\newpage
Non-boundary conforming methods offer an alternative approach to geometry representation.
In this paper we will consider two methods, although there are even more established approaches.
The fictitious boundary method~(FBM)~\cite{Turek2003} is an implicit non-boundary conforming method that employs an indicator function $\tilde{\alpha}(X)$ to mask points in the fluid domain.

\begin{equation}
	\tilde{\alpha}(X) = \left\{
	\begin{matrix}{}
		1 & \text{for } X \in \Omega_f \\
		0 & \text{for } X \in \Omega\backslash\Omega_f
	\end{matrix}\right.
\end{equation}

The Navier-Stokes equations are then solved on the filtered domain $\Omega_f$.

The immersed boundary method~(IBM)~\cite{Khadra2000, Sasse2025} is an explicit non-boundary conforming method that locally blocks the flow by assigning a pseudo-porosity $\alpha(X)$ to each cell of the computational mesh.

\begin{equation}
	\alpha(X) = \left\{
	\begin{matrix}{}
		0 & \text{for } X \in \Omega_f \\
		\alpha_{max} & \text{for } X \in \Omega_s
	\end{matrix}\right.
\end{equation}

The numerical solution is then obtained by solving the extended Navier-Stokes equations on $\Omega = \Omega_f \cup \Omega_s$:

\begin{align}
	\nabla \cdot \mathbf{u} &= 0 \\
	\rho \left(\mathbf{u} \cdot \nabla \right) \mathbf{u} - \nabla \cdot \sigma(\mathbf{u,p}) \ + \alpha \mathbf{u} &= 0
\end{align}

Both non-boundary conforming methods enable computationally cheap mesh generation that is independent of the geometric complexity of the flow channel. However, this comes at a loss of accuracy, particularly at the interface.
While the interface in FBM can be recovered through evaluation of the gradient of the indicator function $\nabla\tilde{\alpha}(X)$, that information is lost in IBM.
This, however, is not necessarily an advantage for FBM, since any data projected to the interface are prone to artifacts in areas with insufficient co-planarity between the background mesh and the interface.
As a result, IBM can sometimes be more robust regarding interface artifacts or instabilities. However, this comes at the additional cost of interface reconstruction for later processing of the optimization result.

\subsection{Adjoint-based Optimization}\label{sec:Opt}

The goal of adjoint-based optimization is to minimize an objective functional $J$ under the PDE constraints $R$, which in this case arise from the incompressible steady Navier-Stokes equations:

\begin{equation}
	\min_\mathbf{b} J(\chi, \mathbf{b}) \quad \text{s. t. } \mathbf{R}(\chi, \mathbf{b}) = \mathbf{0}
\end{equation}

The idea is to introduce a Lagrange functional $L$ with Lagrange multipliers~$\lambda$, also called the adjoint variables:

\begin{equation}
	L := J + \int_{\Omega} \lambda R d\Omega
\end{equation}

This allows for an optimization of a suitable design variable $\mathbf{b}$ through a gradient descent method (Eq.~\ref{eq:gradientdescent}), where the gradient is provided by the directional derivative of the Lagrange functional with respect to the design variable, also known as the design sensitivity:

\begin{equation}
	\mathbf{b}_{n+1} = \mathbf{b}_n - \Delta \frac{\partial L(\mathbf{b}_n)}{\partial \mathbf{b}_n} \label{eq:gradientdescent}
\end{equation}

This directional derivative of the Lagrange functional can be calculated after solving the primal state equations and their adjoint problem in only two solver calls~\cite{Bletsos2023,Othmer2008,Sasse2025}.
The primal problem can be stated as: Find $(\mathbf{u}, p, \mu)$ of sufficient regularity holding

\begin{align}
	\rho(\mathbf{u}\cdot \nabla) \mathbf{u} - \nabla \cdot \sigma(\mathbf{u}, p, \mu) \textcolor{red}{+\alpha \mathbf{u}} &= 0, & \text{in } \; \Omega, \label{eq:primalMomentum} \\
	\nabla \cdot \mathbf{u} &= 0 & \text{in} \; \Omega, \\
	\mu - F(\dot{\gamma}) &= 0 & \text{in} \; \Omega, \label{eq:primalVisc} \\
	\mathbf{u} &= \mathbf{u}_0 & \text{on} \; \partial\Omega_{\text{in}}, \\
	\mathbf{u} &= 0 & \text{on} \; \partial\Omega\backslash(\partial\Omega_{\text{in}}\cup \partial\Omega_{\text{out}}), \\
	(2\mu\mathbf{D}(\mathbf{u})-p\mathbf{I})\cdot \mathbf{n} &= 0 & \text{on} \; \partial\Omega_{\text{out}}, \\
	\sigma(\mathbf{u}, p, \mu) &:= (2\mu \mathbf{D}(\mathbf{u}) - \mathbf{I}p), &
\end{align}

where Eq.~\ref{eq:primalVisc} is an auxiliary problem needed for the shear-rate dependent viscosity in the adjoint model~\cite{Bletsos2023}. 
\newpage
The corresponding adjoint model reads:
Given the primal solution $(\mathbf{u}, p, \mu)$, find $(\mathbf{u}_a, p_a)$ holding

\begin{align}
	&\frac{\partial j_\Omega}{\partial \mathbf{u}} + \rho \left(\mathbf{u} \cdot \nabla\right)\mathbf{u_a} - \rho \left(\mathbf{u_a} \cdot \nabla\right)\mathbf{u} -\sigma(\mathbf{u_a}, p_a, \mu) + \nabla \cdot \left(\mu_a \mathbf{X}\right) \textcolor{red}{+\alpha \mathbf{u}} =0  &\text{in} \; \Omega, \label{eq:adjointMomentum} \\
	&\nabla \cdot \mathbf{u_a} - \frac{\partial j_{\Omega}}{\partial p} = 0  &\text{in} \; \Omega, \\
	&\frac{\partial j_{\Omega}}{\partial \mu} + 2\mathbf{D}(\mathbf{u}) :\mathbf{D}(\mathbf{u_a}) +\mu_a  = 0 & \text{in} \; \Omega, \\
	&\mathbf{u_a}|_{t} = 0, \; \mathbf{u_a}|_{\mathbf{n}} = -\frac{\partial j_{\Gamma}}{\partial p}  &\text{on} \; \partial\Omega\backslash\partial\Omega_{\text{out}}, \\
	&\sigma(\mathbf{u_a}, p_a, \mu)\cdot \mathbf{n} = - \frac{\partial j_{\Gamma}}{\partial \mathbf{u}} - \rho\mathbf{u_a}(\mathbf{u} \cdot \mathbf{n}) + \mu_a \mathbf{X} \cdot \mathbf{n}  &\text{on} \; \partial\Omega_{\text{out}}, \\
	&\mu = F(\dot{\gamma}), & \\
	&\mathbf{X} := 2\dot{\gamma}^{-1}F'(\dot{\gamma})\mathbf{D}(\mathbf{u}), & \\
	&\mu_a = - 2\mathbf{D}(\mathbf{u}):\mathbf{D}(\mathbf{u_a}) - \frac{\partial j_\Omega}{\partial \mu}. & 		
\end{align}

Note that the highlighted term $+ \alpha \mathbf{u}$ in Eq.~\ref{eq:primalMomentum} and Eq.~\ref{eq:adjointMomentum} is only added when using IBM.
The solution of the primal and adjoint system is then used to calculate the design sensitivity for the optimization of the flow channel geometry.

Adjoint-based optimization can be performed using two approaches, shape optimization and topology optimization~\cite{Othmer2008}.
In shape optimization, the design variable is the localized surface normal displacement $\beta$ on the nodes of the design boundary. The corresponding surface sensitivity $\sigma_{\beta}$ is calculated as

\begin{equation}
	\delta_b L \cdot \delta \mathbf{b} = \sigma_{\beta} = \left[-2\mu\mathbf{D}(\mathbf{u_a}) \cdot\mathbf{n} + \mu_a 2\dot{\gamma}^{-1}F'(\dot{\gamma})\mathbf{D}(\mathbf{u})\cdot \mathbf{n}  \right]\cdot \frac{\partial \mathbf{u}}{\partial \mathbf{n}}.
\end{equation}

The shape update is then performed by moving the points on the design boundary according to the localized surface normal displacements $\beta$ and performing a subsequent smoothing step on the inner nodes of the mesh.
Topology optimization, on the other hand, offers more design freedom since it is not restricted to the initial topology of the design boundary. For our purpose, the design variable is the IBM variable $\alpha$ and the corresponding volume sensitivity $\sigma_{\alpha_i}$ for each cell $i$ of the computational mesh is calculated using the local cell volume $V_i$ with

\begin{equation}
	\delta_b L \cdot \delta \mathbf{b} = \sigma_{\alpha_i} = u_{a,i} \cdot u_i V_i. \label{eq:VolSenseCalc}
\end{equation}

The subsequent topology update ($n \rightarrow n+1$) is performed using a gradient-based update scheme

\begin{equation}
	\alpha_{n+1} = \alpha_{n} (1+\gamma) + \gamma \min (\max (\alpha_{n} - \lambda^* \sigma_{\alpha}, 0), \alpha_{max}) \label{eq:VolSenseUpdate}
\end{equation}

where $\gamma$ and $\lambda^*$ refer to an underrelaxation factor and the response step size, respectively.
While this approach offers more design freedom, it can also require additional geometric constraints to ensure optimization results that are also viable from an engineering perspective.

\newpage
\subsection{Discussion}

The geometry representation strategies presented in Sec.~\ref{sec:Geometry} and the optimization algorithms presented in Sec.~\ref{sec:Opt} allow for different algorithm pairings. However, not all of them are equally suitable for a successful optimization for profile extrusion dies.
Both shape optimization using an IBM approach and topology optimization using an FBM approach are not feasible within the described framework since the involved components are inherently incompatible.

\begin{figure}[htb]
	\centering
	\begin{subfigure}{0.49\textwidth}
		\centering
		\includegraphics[clip, trim= 275 0 350 0, width=0.99\textwidth]{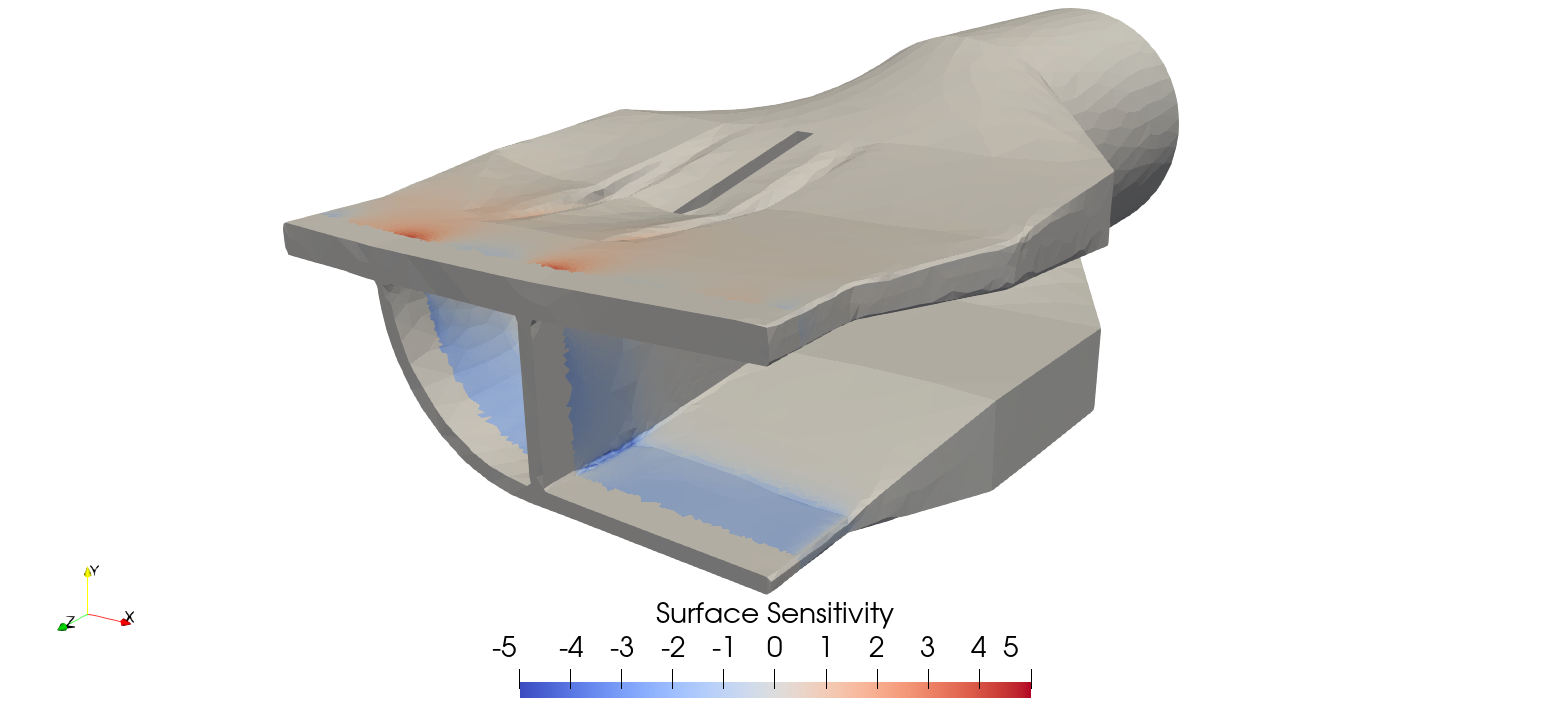}
		\caption{Surface sensitivities calculated with aligned mesh.}
		\label{fig:AlignedSurfaceSense}
	\end{subfigure}
	\hfill
	\begin{subfigure}{0.49\textwidth}
		\centering
		\includegraphics[clip, trim= 275 0 350 0, width=0.99\textwidth]{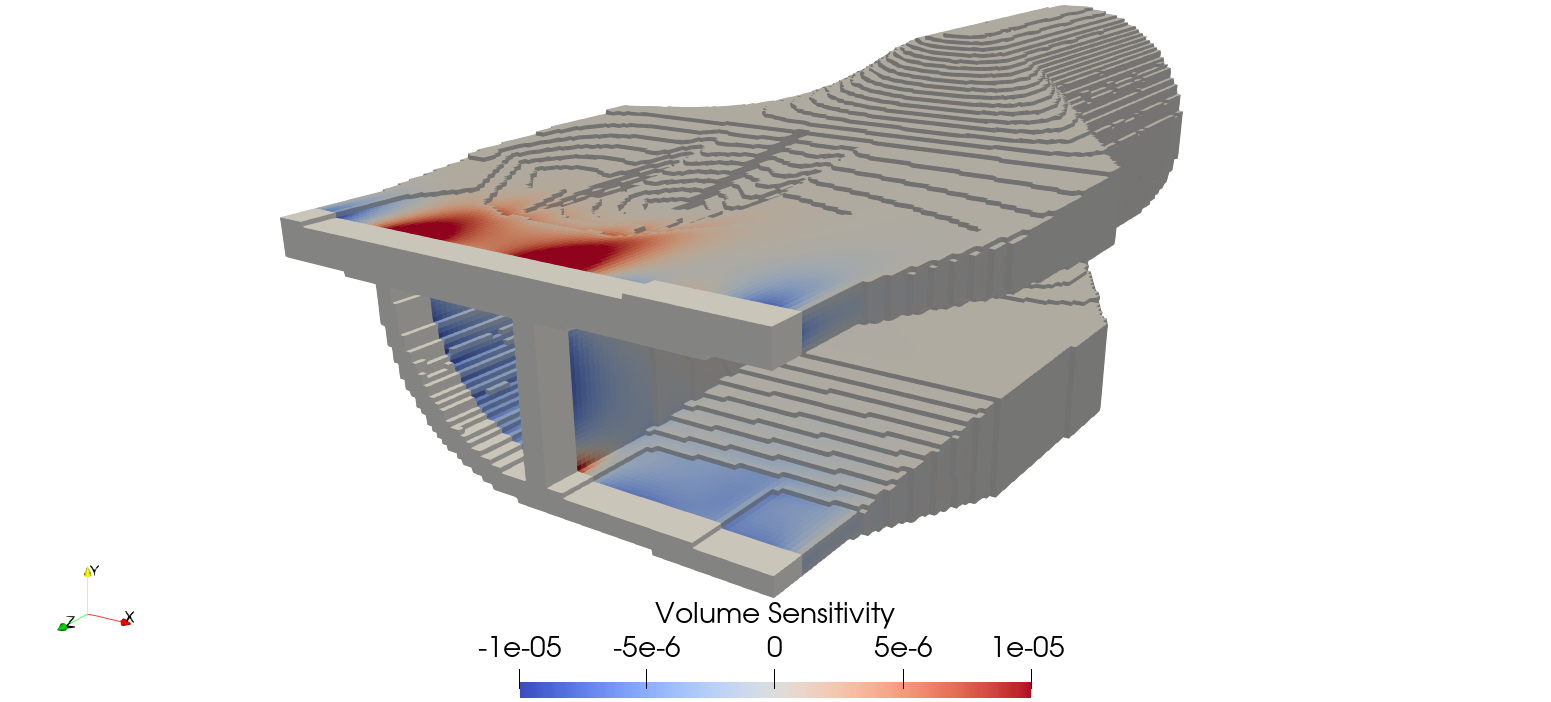}
		\caption{Volume sensitivities calculated using IBM.}
		\label{fig:IBMVolumeSense}
	\end{subfigure}
	\centering
	\begin{subfigure}{0.49\textwidth}
		\centering
		\includegraphics[clip, trim= 275 0 350 0, width=0.99\textwidth]{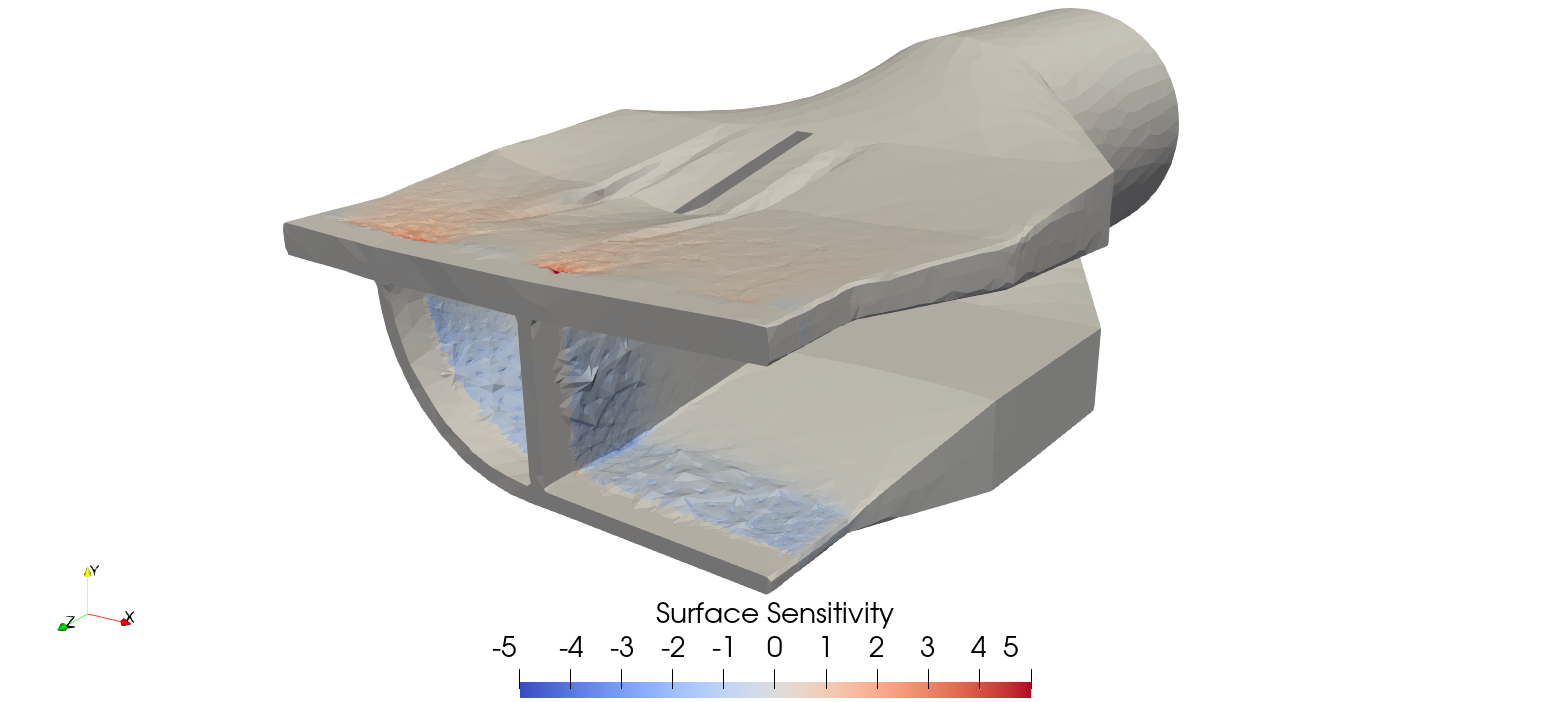}
		\caption{Resulting flow channel update with aligned mesh.}
		\label{fig:AlignedSurfaceSenseUpdate}
	\end{subfigure}
	\hfill
	\begin{subfigure}{0.49\textwidth}
		\centering
		\includegraphics[clip, trim= 275 0 350 0, width=0.99\textwidth]{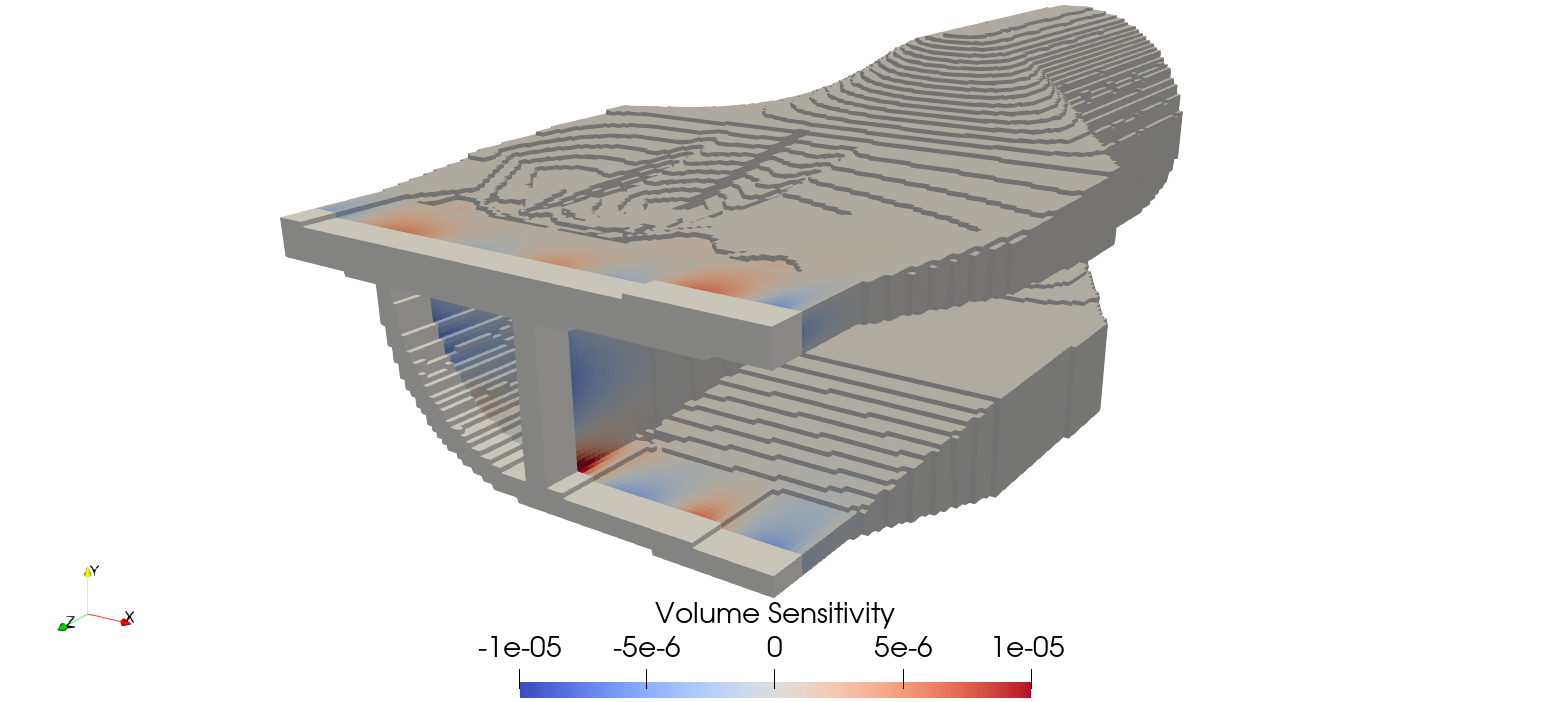}
		\caption{Resulting flow channel update using IBM.}
		\label{fig:IBMVolumeSenseUpdate}
	\end{subfigure}
	\caption{Suitable geometry modeling and optimization strategy pairings for a complex flow channel.}
	\label{fig:SuitablePairingsProfileKTP}
\end{figure}

The most intuitive methods are shape optimization using an aligned mesh (see Fig.~\ref{fig:AlignedSurfaceSense}) and topology optimization using IBM (see Fig.~\ref{fig:IBMVolumeSense}).
While the former leads to a visually superior sensitivity distribution, the optimization result will depend a lot on the computational mesh.
The difference becomes apparent in Fig.~\ref{fig:AlignedSurfaceSenseUpdate} and Fig.~\ref{fig:IBMVolumeSenseUpdate}, where the pattern of the underlying computational mesh dictates the shape of the optimized surface.
As a result, any aligned mesh used for this purpose has to adhere to even more restrictions than are usually imposed on meshes for computational fluid dynamics.
This makes topology optimization using IBM the less "pretty" but more robust choice for an optimization algorithm for complex geometries.

\begin{figure}[htb]
	\centering
	\begin{subfigure}[b]{0.49\textwidth}
		\centering
		\includegraphics[clip, trim= 300 126 350 40, width=0.99\textwidth]{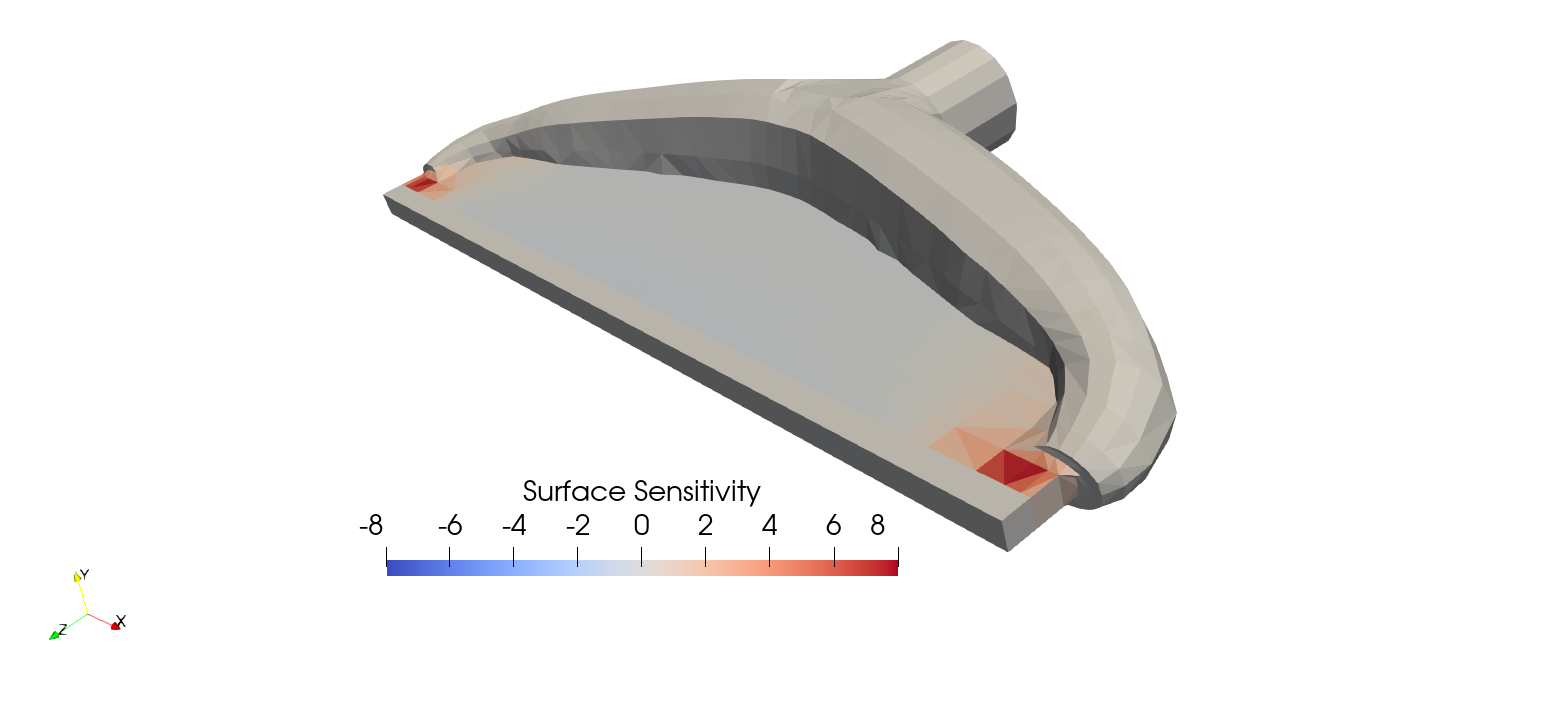}
		\caption{Surface sensitivities on an aligned mesh.\newline \textcolor{white}{.}}
		\label{fig:FlatsheetAlignedSurfaceSense}
	\end{subfigure}
	\hfill
	\begin{subfigure}[b]{0.49\textwidth}
		\centering
		\includegraphics[clip, trim= 300 126 350 40, width=0.99\textwidth]{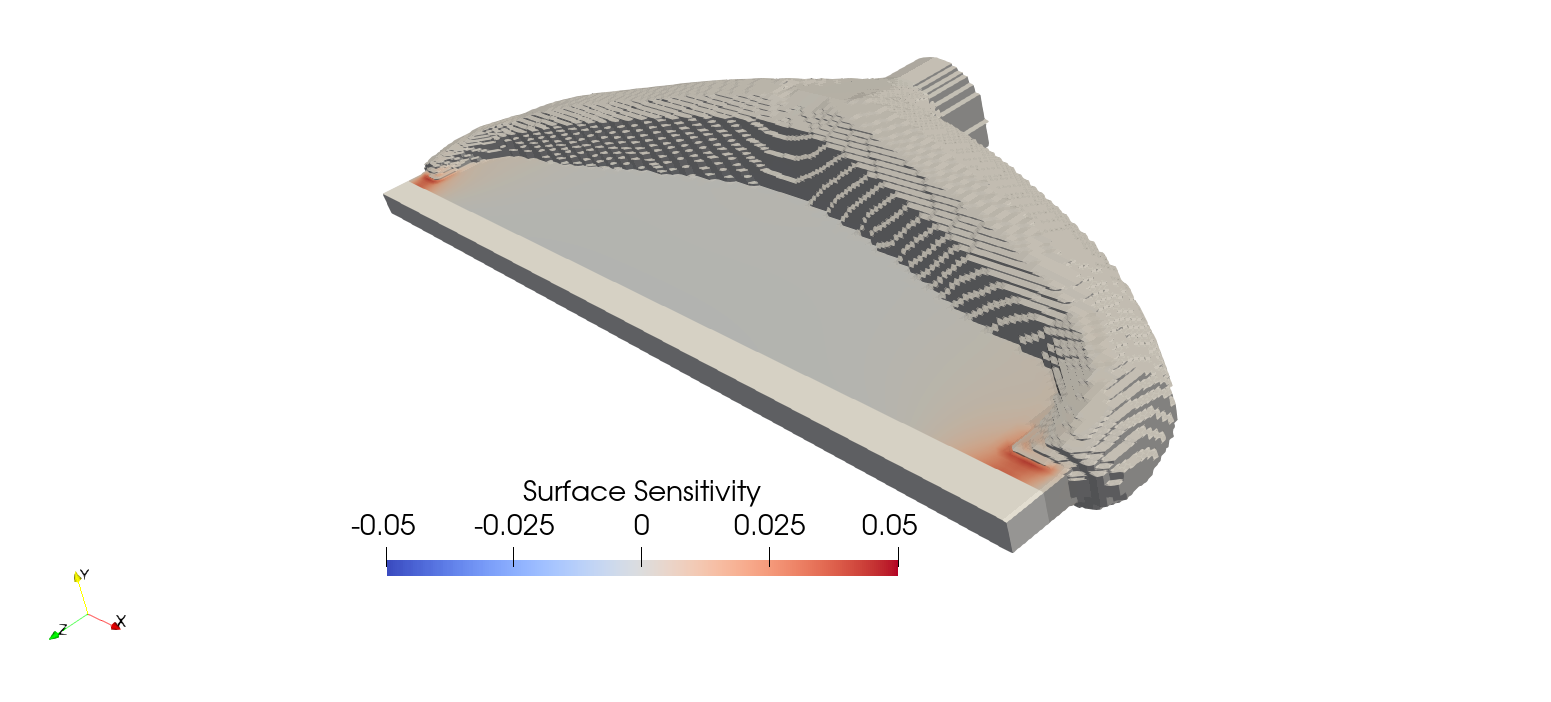}
		\caption{Surface sensitivities calculated using FBM\newline with a co-planar interface.}
		\label{fig:FlatsheetFBMSurfaceSense}
	\end{subfigure}
	\centering
	\begin{subfigure}[b]{0.49\textwidth}
		\centering
		\includegraphics[clip, trim= 275 0 325 0, width=0.99\textwidth]{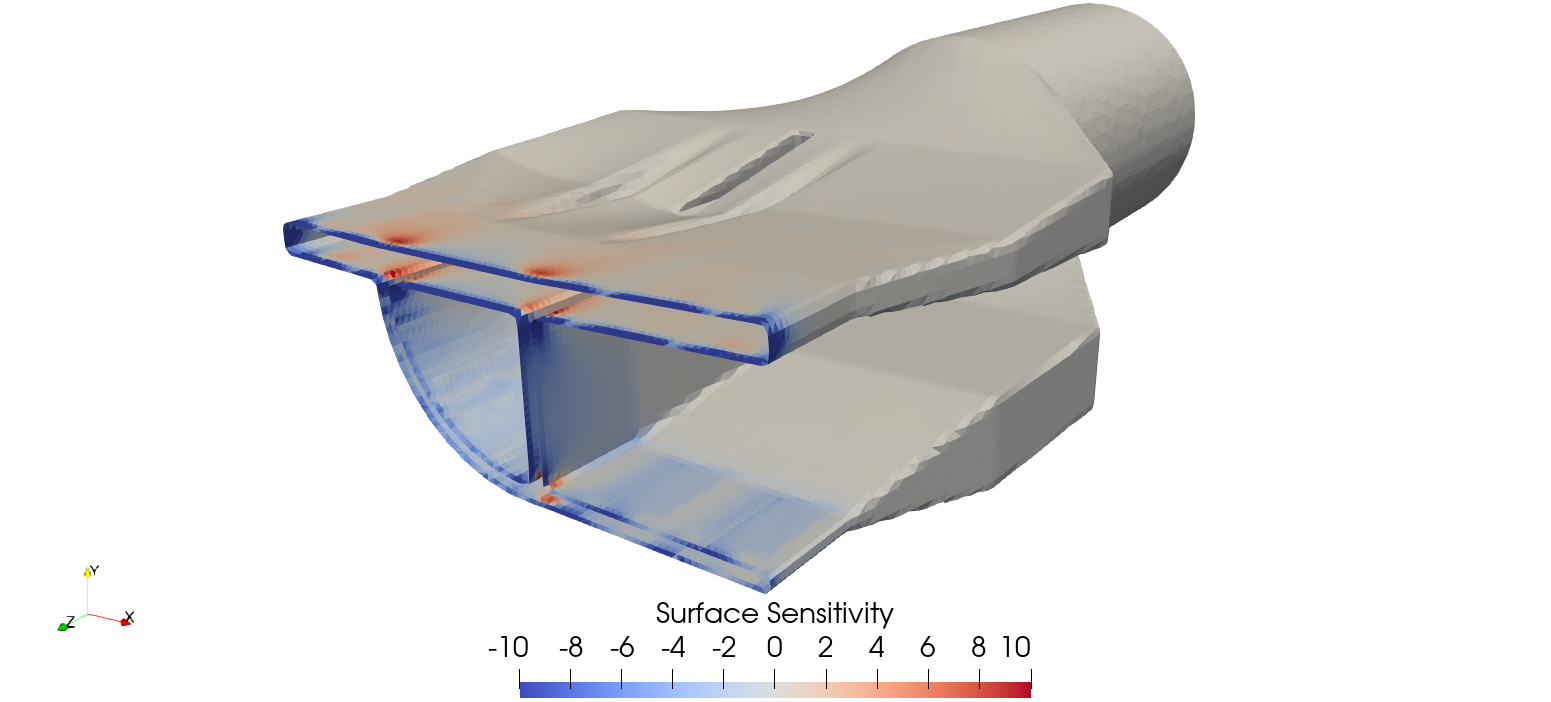}
		\caption{Surface sensitivities calculated using FBM\newline with no co-planar interface.}
		\label{fig:FBMSurfaceSense}
	\end{subfigure}
	\hfill
	\begin{subfigure}[b]{0.49\textwidth}
		\centering
		\includegraphics[clip, trim= 275 0 325 0, width=0.99\textwidth]{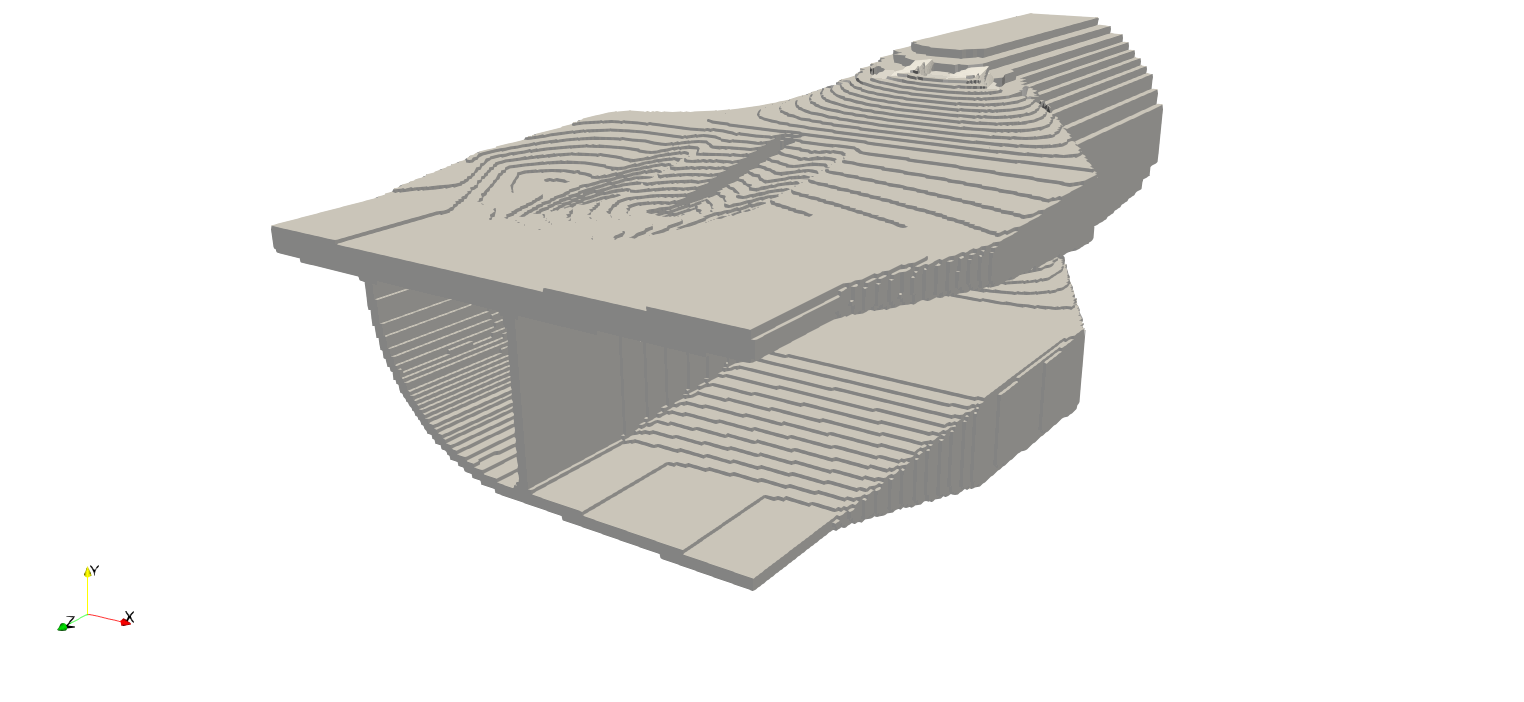}
		\caption{Staircase pattern in filtered FBM domain\newline due to no co-planar interface.}
		\label{fig:FBMStaircase}
	\end{subfigure}
	\caption{Example where FBM leads to artifacts in the calculated surface sensitivity.}
	\label{fig:NotSuitablePairingsProfileKTP}
\end{figure}

Topology optimization with an aligned mesh can work. However, this requires an appropriate interface resolution and is usually not feasible for complex geometries.
Shape optimization using an FBM geometry representation can be successful.
Fig.~\ref{fig:FlatsheetAlignedSurfaceSense} and Fig.~\ref{fig:FlatsheetFBMSurfaceSense} depict the surface sensitivity for a flatsheet extrusion die using both an aligned mesh and FBM.
However, this approach breaks down as soon as the interface is not co-planar with the background mesh.
Revisiting the example given in Figs.~\ref{fig:AlignedSurfaceSense}-\ref{fig:IBMVolumeSense}, the surface of the flow channel is not co-planar with the background mesh (see Fig.~\ref{fig:FBMStaircase}). As a result, the projection of the local velocity to the interface leads to wave-patterned artifacts (see Fig.~\ref{fig:FBMSurfaceSense}).
While these can be reduced by smoothing, this approach risks removing real results in the process.

As a result, the immersed boundary method combined with a topology optimization algorithm appears to be the most versatile method for use cases with industrial relevance.
However, shape optimization with aligned meshes will still produce the most reliable result, provided an appropriate aligned mesh is available.


\section{Application to Profile Extrusion}

\subsection{Setup and Numerical Model}

\begin{figure}[htb]
	\centering
	\begin{subfigure}{0.54\textwidth}
		\centering
		\includegraphics[clip, trim = 200 105 150 55, width=\textwidth]{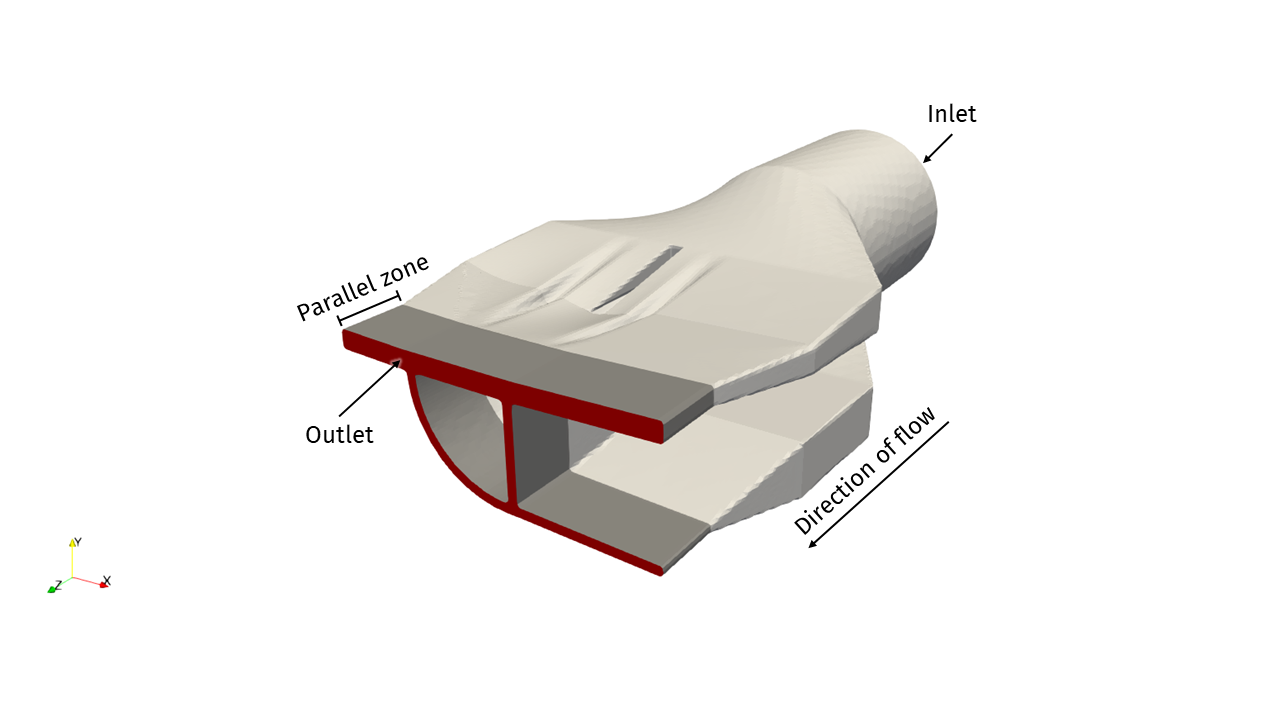}
		\caption{Initial flow channel geometry.}
		\label{fig:ProfileKTPSetup}
	\end{subfigure}
	\hfill
	\begin{subfigure}{0.44\textwidth}
		\centering
		\includegraphics[clip, trim = 125 90 190 70, width=\textwidth]{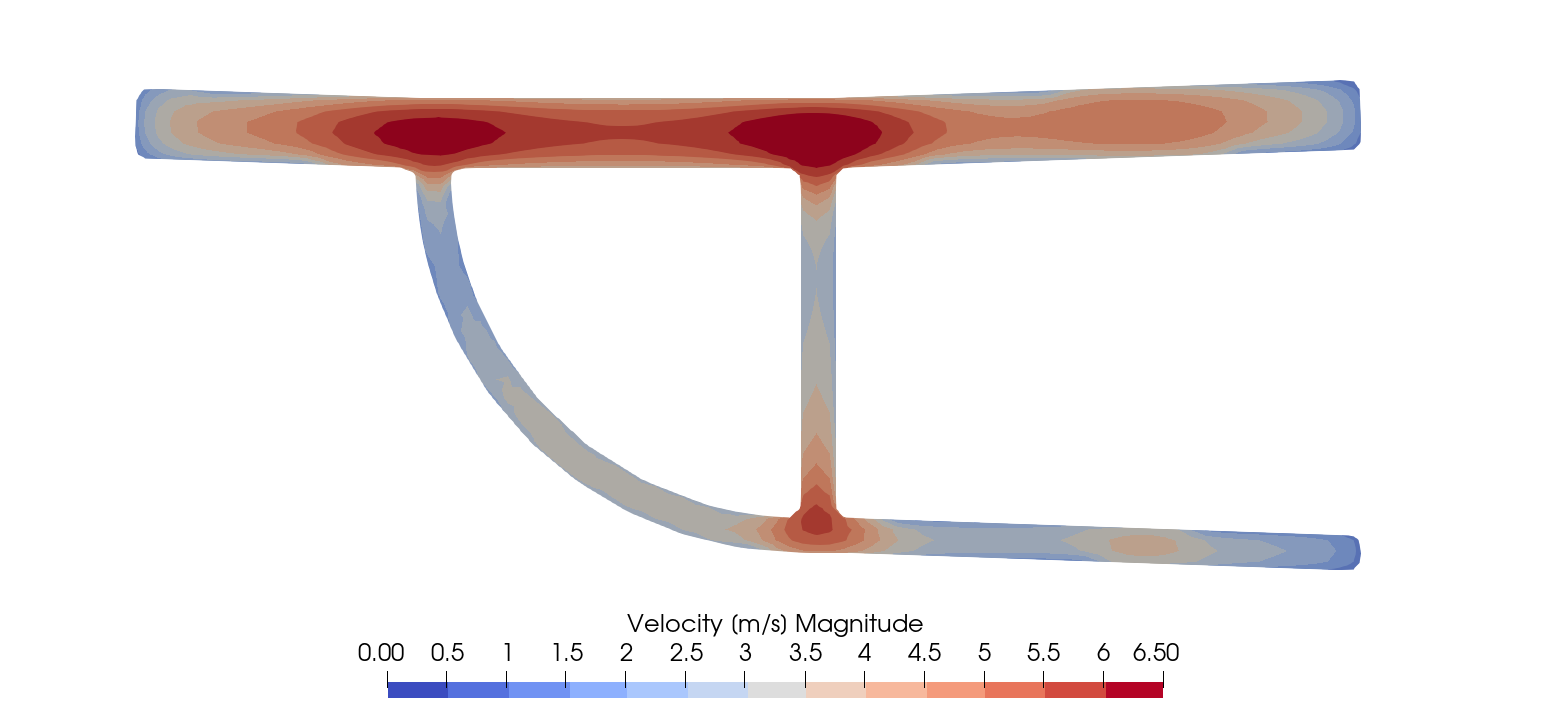}
		\caption{Velocity distribution at outlet before optimization.}
		\label{fig:ProfileKTPOutflow}
	\end{subfigure}
	\caption{Setup of the demonstrator profile extrusion die.}
	\label{fig:DemonstratorCaseSetup}
\end{figure}

This flow channel geometry is based on a real profile extrusion die for an extruder with a 45~mm diameter. It shapes the melt into the complex outflow profile depicted in Fig.~\ref{fig:ProfileKTPSetup}.
In addition to other typical features in traditional extrusion die design, it includes a 15~mm parallel zone near the outlet where the geometry cannot be modified in order to control the die swell by allowing for sufficient melt relaxation.
The initial and boundary conditions for this use case are also taken from a real process operating point with a throughput of 10~kg/h and a processing temperature of 220~°C.
The extruded material is a polyethylene of high density (PE-HD) and is modeled with a Carreau-WLF approach.
The first optimization objective for this use case is the maximization of the flow balance at the outlet (initial configuration shown in Fig.~\ref{fig:ProfileKTPOutflow}), which is equivalent to the minimization of the deviation of the local velocity from the average velocity at the outlet $\mathbf{\bar{u}}$.

\begin{equation}
	\min\ J = \int_{\Gamma_{out}} \frac{1}{2} \left(\mathbf{u} - \mathbf{\bar{u}} \right)^2
\end{equation}

\begin{figure}[htb]
	\centering
	\includegraphics[clip, trim= 270 205 150 0, width=0.44\textwidth]{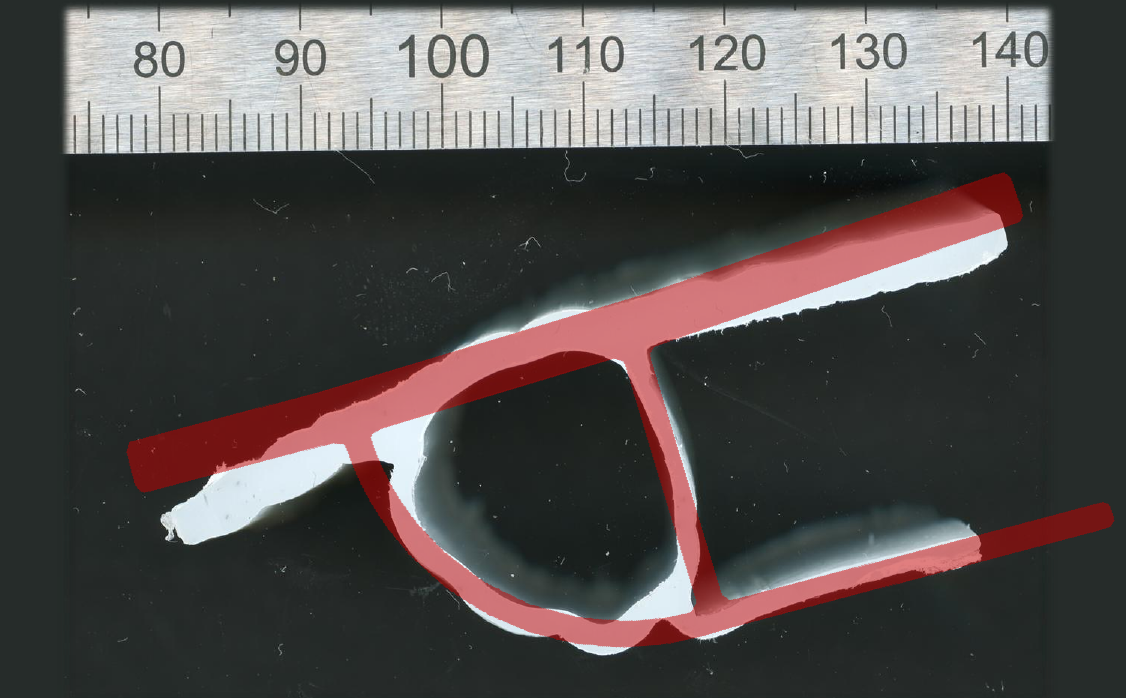}
	\caption{Difference between the target profile (red) and the actual profile (white) in extrusion trials before flow channel optimization. Image provided by Kunststofftechnik Paderborn~(KTP).}
	\label{fig:realFlowChannel}
\end{figure}

The imbalance in the outflow velocity distribution is exacerbated by the high viscosity and shear-thinning behavior of the extruded material. In practice, this means that the flow channel of the die is not filled appropriately (see Fig.~\ref{fig:realFlowChannel}). 

The flow channel geometry was modeled using an IBM approach, embedding the geometry in a cartesian background mesh. $\alpha_{max}$ was chosen at $10^5$, balancing blockage of the flow with numerical stability.
The optimization algorithm calculates the volume sensitivity according to Eq.~\ref{eq:VolSenseCalc}.
The topology update is effectively reduced to a shape update by restricting the geometry update in Eq.~\ref{eq:VolSenseUpdate} to cells at the interface, i.~e., cells that have vertices on the fluid and the solid side of the interface.
No underrelaxation was used in this case ($\gamma = 1$) and the response step size was manually fitted to $\lambda^* = 10^{10}$. A total of 20 optimization loops were performed, meaning that the results shown later will represent the optimization result after 19 optimization steps.
The resulting flow channel geometry is reconstructed using a SurfaceNets algorithm~\cite{Gibson1998}, ensuring a 2-manifold surface using the approach by \textit{Schaefer et al.}~\cite{Schaefer2007}, and exported as a surface mesh file.

The numerical model was implemented in the open-source finite element software FEAT3, which is part of the FEATFLOW software family~\cite{FEAT3}.
The spatial discretization is performed using $Q^2$/$P_{disc}^1$ elements, which are inf-sup stable and particularly suited for non-Newtonian flows without the need for additional stabilization~\cite{Boffi2013}.
The numerical model is solved as follows~\cite{Esser2026}: 
The nonlinear primal system is solved directly into its steady state using a monolithic Newton solver. The adjoint system does not require this step since it is already linear.
The linearized systems are then solved using a monolithic geometric multigrid solver.
Using a background mesh with 3,750 hexahedral elements on the coarsest level with three refinement levels, this results in a total of $\sim$48,000,000 degrees of freedom per optimization loop.
The solver uses a Vanka-type cell-wise additive Schwarz smoother with a GMRES solver at the coarsest level. 
Simulations were run using two compute servers, each containing two AMD EPYC 9354 CPUs with a total of 64 cores and 1.5 TB of DDR5 @ 4800 MHz RAM, resulting in a run time of 20~minutes in the first optimization loop.
Later optimization steps are initialized using the solutions from the previous step to reduce computational cost.

\subsection{Results}

\begin{figure}[p]
	\centering
	\begin{subfigure}[b]{0.49\textwidth}
		\centering
		\includegraphics[clip, trim = 300 0 350 25, width=\textwidth]{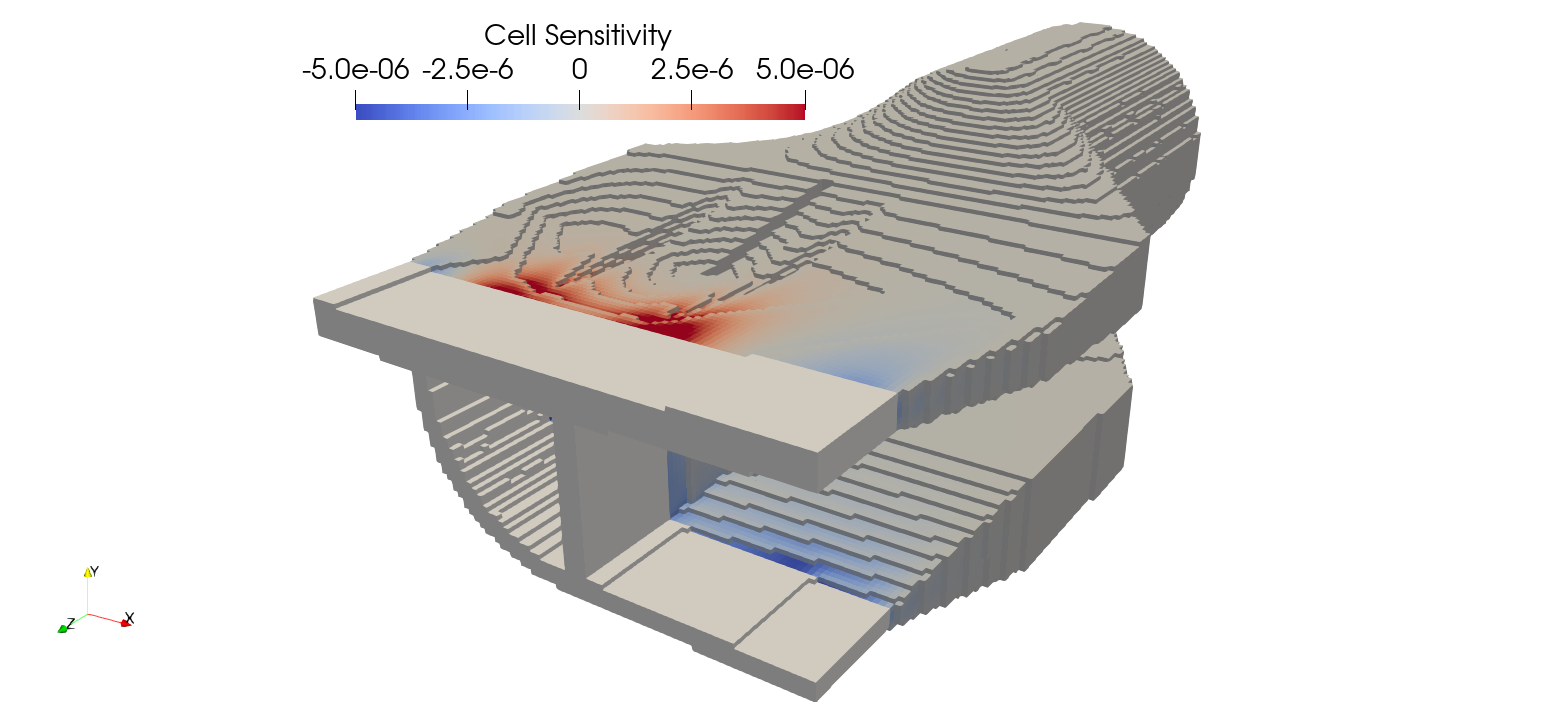}
		\caption{Volume sensitivity before optimization.}
		\label{fig:SensitivityBefore}
	\end{subfigure}
	\hfill
	\begin{subfigure}[b]{0.49\textwidth}
		\centering
		\includegraphics[clip, trim = 300 0 350 25, width=\textwidth]{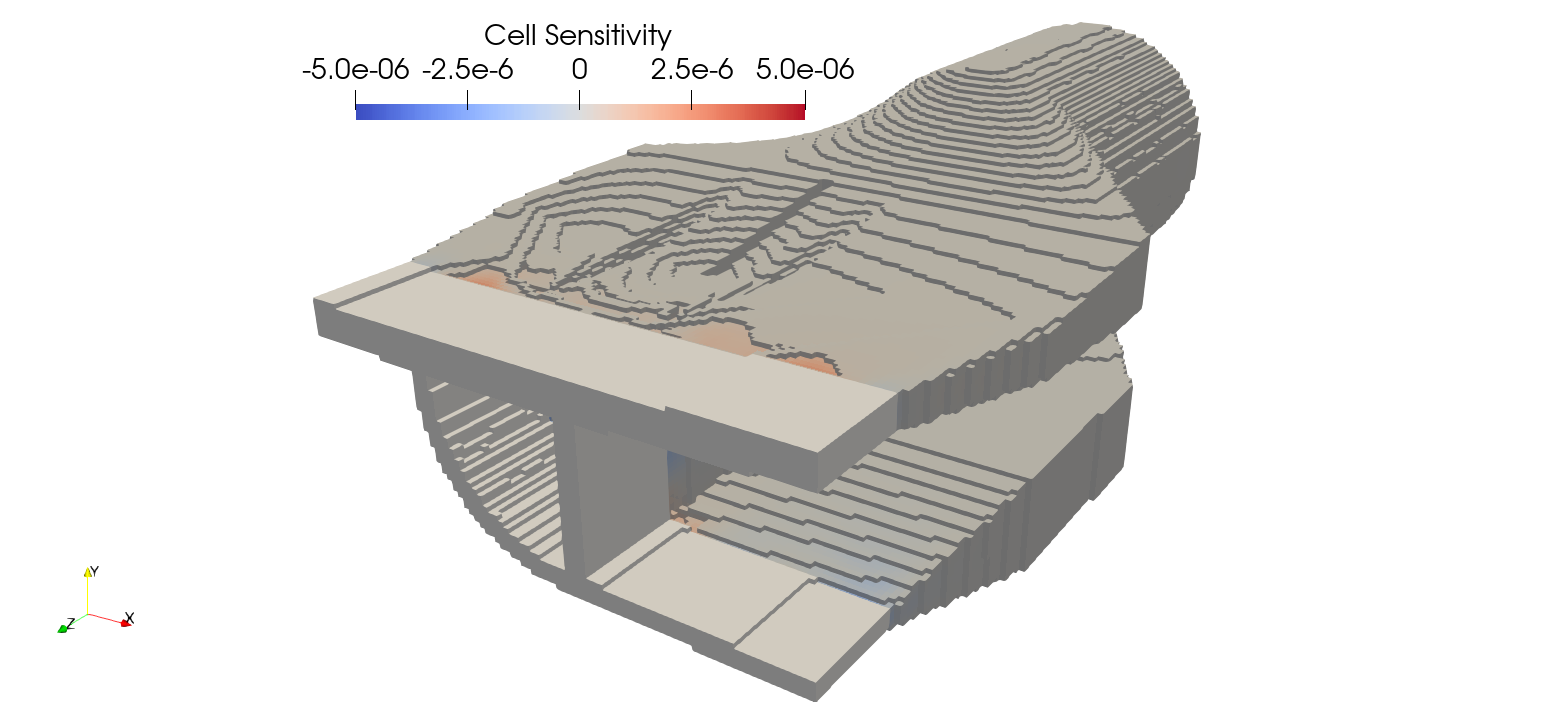}
		\caption{Volume sensitivity after 19 optimization steps.}
		\label{fig:SensitivityAfter}
	\end{subfigure}
	\begin{subfigure}[b]{0.49\textwidth}
		\centering
		\includegraphics[clip, trim= 190 50 200 50, width=\textwidth]{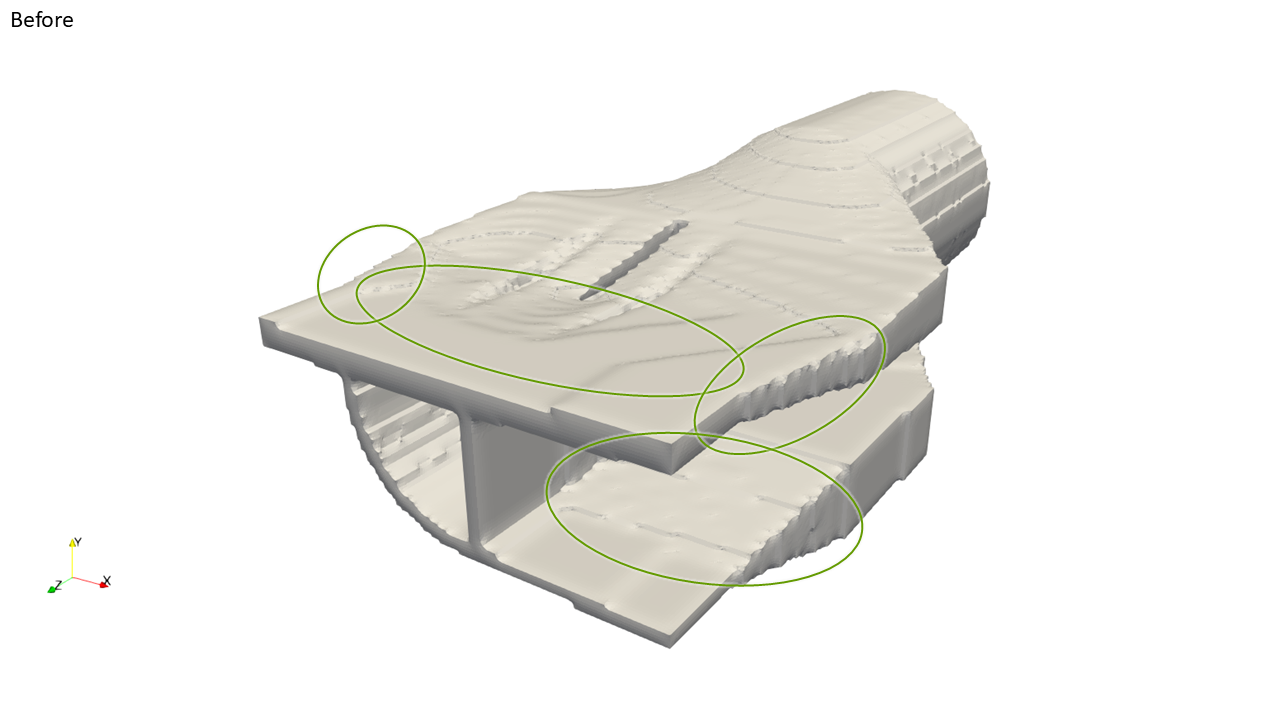}
		\caption{Flow channel geometry before optimization.}
		\label{fig:GeometryBefore}
	\end{subfigure}
	\hfill
	\begin{subfigure}[b]{0.49\textwidth}
		\centering
		\includegraphics[clip, trim= 190 50 200 50, width=\textwidth]{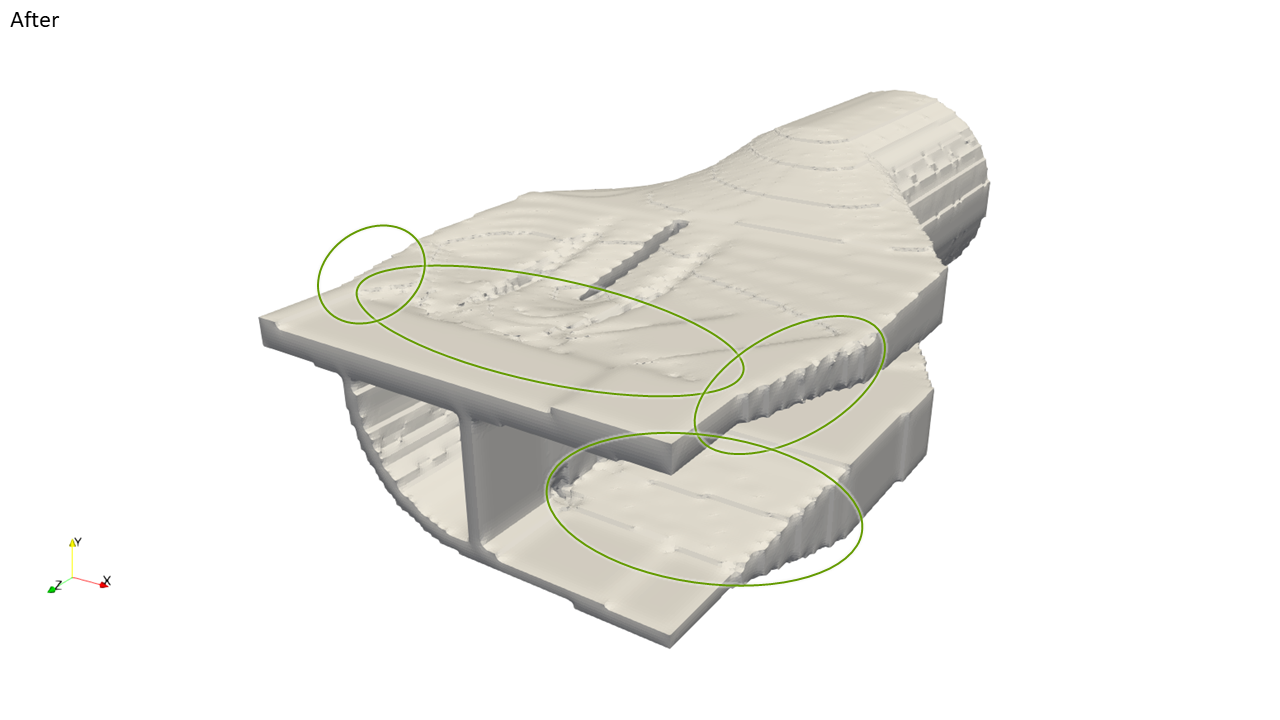}
		\caption{Flow channel geometry after 19 optimization steps.}
		\label{fig:GeometryAfter}
	\end{subfigure}
	\begin{subfigure}[b]{0.49\textwidth}
		\centering
		\includegraphics[clip, trim = 125 0 190 70, width=\textwidth]{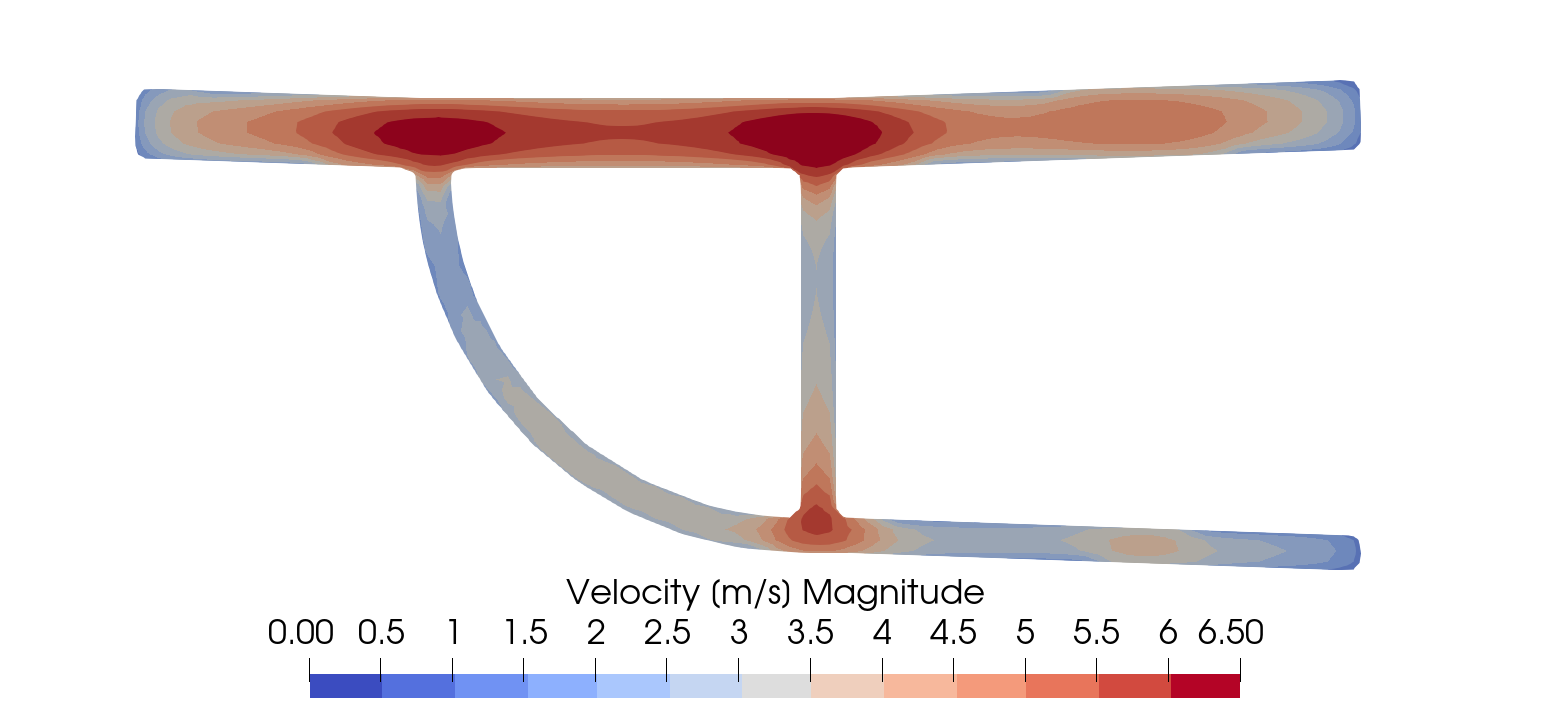}
		\caption{Velocity distribution at outlet before optimization.\newline \textcolor{white}{.}}
		\label{fig:VelocityBefore}
	\end{subfigure}
	\hfill
	\begin{subfigure}[b]{0.49\textwidth}
		\centering
		\includegraphics[clip, trim = 125 0 190 70, width=\textwidth]{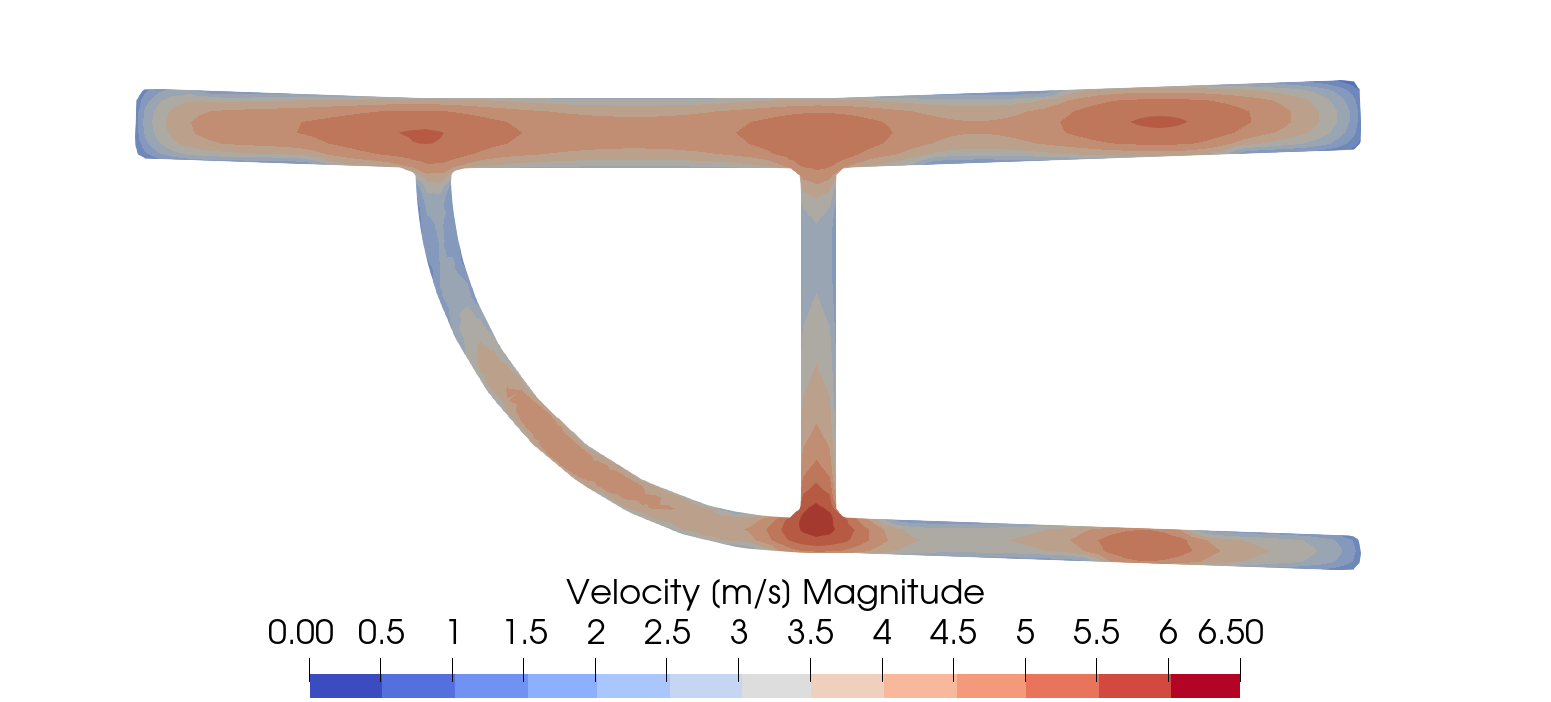}
		\caption{Velocity distribution at outlet after \newline 19 optimization steps.}
		\label{fig:VelocityAfter}
	\end{subfigure}
	\caption{Changes to flow channel geometry and resulting changes to outflow velocity distribution.}
	\label{fig:OptimizationResults}
\end{figure}

Fig.~\ref{fig:SensitivityBefore} depicts the resulting volume sensitivities in the first optimization loop. Positive sensitivities (red areas) highlight areas where the algorithm wants to restrict the flow channel, while negative sensitivities (blue areas) highlight parts of the flow channels that should be expanded to improve the flow balance.
The former primarily affects the flow channel near the upper T-junction, where the outflow velocities are too high, while the latter mainly affects the lower parts of the flow channel, where the small channel widths previously led to an incomplete filling of the die's flow channel (see Fig.~\ref{fig:realFlowChannel}).
Figs.~\ref{fig:GeometryBefore}-\ref{fig:GeometryAfter} show how the flow channel was adapted according to Eq.~\ref{eq:VolSenseUpdate}, reducing the cost functional from 439.435 to 207.800 in the process.
This results in a reduction of the standard deviation of the velocity at the outlet from 2507.75 to 1128.65, which in turn means an improvement in the flow balance by 55\% (see Figs.~\ref{fig:VelocityBefore}-\ref{fig:VelocityAfter}).
Although the optimization algorithm has not fully converged at this point, the most significant changes to the shape are finished; the remaining areas with potential for optimization focus on the parts of the flow channel that could still be expanded a bit more (see Fig.~\ref{fig:SensitivityAfter}).
%

\section{Conclusion and Outlook}

The results of this study demonstrate that adjoint-based optimization can be successful using various pathways, but not all of them are equally suited for engineering applications. In particular, not all methods suited for CFD analysis are automatically suited for shape or topology optimization. For complex profile extrusion dies, topology optimization using an immersed boundary method is the most versatile approach.
Applying the optimization algorithm to a complex profile extrusion die led to an improvement of the flow balancing by 55\% using fewer solver calls than other established CFD-based optimization methods.

Future work will focus on the investigation of further improvements to geometry representation, e.g., through adaptive mesh refinement or local deformation of the mesh \cite{Grajewski2009}. In addition, the developed optimization algorithm will be extended to fully non-isothermal flows and routines improving the robustness of the optimization algorithm will be added.
The final step in this research will be the manufacturing of an optimized profile extrusion die based on these optimization results for subsequent validation in lab trials.

\section*{Acknowledgments}

This project is supported by the Federal Ministry for Economic Affairs and Energy (BMWE) on the basis of a decision by the German Bundestag.

\end{document}